\documentclass[12pt,letterpaper,english]{article}
\usepackage[T1]{fontenc}
\usepackage[latin9]{inputenc}
\usepackage{float}
\usepackage{amsthm}
\usepackage{amsmath}
\usepackage{graphicx}
\usepackage{setspace}
\usepackage{amssymb}
\makeatletter

\floatstyle{ruled}
\newfloat{algorithm}{tbp}{loa}
\floatname{algorithm}{Algorithm}

\newenvironment{lyxcode}
{\par\begin{list}{}{
\setlength{\rightmargin}{\leftmargin}
\setlength{\listparindent}{0pt}% needed for AMS classes
\raggedright
\setlength{\itemsep}{0pt}
\setlength{\parsep}{0pt}
\normalfont\ttfamily}%
 \item[]}
{\end{list}}

\author{Sheng Yu\\ Advisor: Enrique Campos-N\'{a}\~{n}ez\\ \normalsize Department of Engineering Management \& System Engineering\\ \normalsize The George Washington University\\}
\date{}

\makeatother

\usepackage{babel}

\begin{document}

\title{Introducing the Adaptive Convex Enveloping}
\maketitle
\begin{abstract}
Convexity, though extremely important in mathematical programming,
has not drawn enough attention in the field of dynamic programming.
This paper gives conditions for verifying convexity of the cost-to-go
functions, and introduces an accurate, fast and reliable algorithm
for solving convex dynamic programs with multivariate continuous states
and actions, called Adaptive Convex Enveloping. This is a short introduction
of the core technique created and used in my dissertation, so it is
less formal, and misses some parts, such as literature review and
reference, compared to a full journal paper.
\end{abstract}

\subsection*{Background Story}

One of my friends has a small firm that designs electric vehicle battery
stations. A battery station works like a gas station: customers arrive
to replace their used batteries with full ones and leave, and it is
the battery station's job to charge the batteries in order to serve
new customers. However, the electricity price fluctuates during the
day. The battery station certainly wants to charge the batteries when
the price is low, but it still has to satisfy customer demand. I wanted
to help my friend develop an optimal charging policy, which is a typical
dynamic program, with multivariate continuous states and actions.
However, current ADP methods are not very satisfying when solving
this class of problems. For example, algorithms generally require
finding the optimal policy under the current approximations of the
cost-to-go functions, but the approximations are in general not convex,
sometimes very wavy, which makes finding the optimal policy difficult
in multivariate problems. Also, even a method is proved to converge
under certain conditions, one generally doesn't know how the current
policy compares to the true optimal one. To better solve this class
of problems, I devised this new method --- Adaptive Convex Enveloping
(A.C.E.).

\subsection*{Formulation}

Let the time periods be $t=1,\dots,T$. Let $x_{t}$, $u_{t}$ and
$w_{t}$ be vectors that represent the state, the action and the random
information at stage $t$, respectively. Let $J_{t}(x_{t})$, $x_{t}\in X_{t}$,
be the cost-to-go function: \begin{align}
J_{t}(x_{t})= & \inf_{u_{t}}E[c_{t}(x_{t},u_{t},w_{t})+J_{t+1}(x_{t+1}(x_{t},u_{t},w_{t}))]\nonumber \\
\text{s.t.}\quad & g_{t}(x_{t},u_{t})\le0,\label{eq:Bellman}\end{align}
where $c_{t}(x_{t},u_{t},w_{t})$ is a known function of the cost
in period $t$, $g_{t}=(g_{1,t},\dots,g_{I_{t},t})^{T}$ is a set
of constraints on the action $u_{t}$, given the current state $x_{t}$,
and $x_{t+1}(x_{t},u_{t},w_{t})$ is the state transition function.

\subsection*{Assumptions}
\begin{itemize}
\item $J_{T}(x_{T})$ is known and convex;
\item The state domain $X_{t}$, $t=1,\dots,T$, is closed, bounded and
convex;
\item The cost function $c_{t}$ with $w_{t}$ fixed, and the constraints
$g_{i,t}$, $i=1,\dots,I_{t}$, are jointly convex on $X_{t}\times\mathbb{R}^{dim(u_{t})}$;
\item If the constraints are not all linear, then the set $\{(x_{t},u_{t})\mid g_{t}(x_{t},u_{t})\le0\}$
must have a relative interior point;
\item The transition $x_{t+1}(x_{t},u_{t},w_{t})$ is linear, i.e., $x_{t+1}=Ax_{t}+Bu_{t}+w_{t}$,
where $A$ and $B$ are known matrices; nonlinear transitions are
not investigated;
\item The distribution of $w_{t}$ is known, and is independent with $x_{t}$
and $u_{t}$.
\end{itemize}

\subsection*{Convexity of the Cost-to-go Functions}

A convex approximation is only desirable when the target function
is convex. So first we have a look at the cost-to-go function. Indeed,
the convexity of the cost-to-go functions is a direct extension of
the results from Fiacco and Kyparisis' paper \textit{Convexity and
Concavity Properties of the Optimal Value Function in Parametric Nonlinear
Programming (1986)}. Note that our action $u_{t}$ is their decision
variable $x$, our state $x_{t}$ is their parameter $\epsilon$,
and our cost-to-go function $J_{t}(x_{t})$ is their optimal value
function $f^{*}(\epsilon)$.

If $J_{t+1}(x_{t+1})$ is convex, then for any realization of $w_{t}$,
$c_{t}(x_{t},u_{t},w_{t})+J_{t+1}(x_{t+1}(x_{t},u_{t},w_{t}))$ is
jointly convex on $X_{t}\times\mathbb{R}^{dim(u_{t})}$ by our assumptions
on $c_{t}$ and the transition, and so is $E[c_{t}(x_{t},u_{t},w_{t})+J_{t+1}(x_{t+1}(x_{t},u_{t},w_{t}))]$.
Since $X_{t}$ is a convex set, and constraints $g_{i,t}$, $i=1,\dots,I_{t}$,
are jointly convex on $X_{t}\times\mathbb{R}^{dim(u_{t})}$, the point-to-set
map $R(x_{t})=\{u_{t}\mid g_{t}(x_{t},u_{t})\le0\}$ is convex on
$X_{t}$ by Proposition 2.3 from Fiacco and Kyparisis (1986). Then
by their Proposition 2.1, $J_{t}(x_{t})$ is convex on $X_{t}$.

Since $J_{T}(x_{T})$ is convex, the above argument implies that $J_{t}(x_{t})$
is convex for all $t=1,\dots,T$.

\subsection*{Convex Enveloping of the Cost-to-go Functions}

Since the cost-to-go functions are convex, we can use supporting hyperplanes
as their (outer) approximations. Assume that $J_{t+1}(x_{t+1})$ has
been investigated at $\{x_{t+1}^{j}\}_{j=1}^{N}$, where we know the
function value $J_{t+1}(x_{t+1}^{j})$, as well as a subgradient $\nabla J_{t+1}(x_{t+1}^{j})$
(with abuse of the notation). Thus for any point $x_{t+1}$, we know
from convexity that \[
J_{t+1}(x_{t+1})\ge J_{t+1}(x_{t+1}^{j})+\nabla J_{t+1}(x_{t+1}^{j})^{T}(x_{t+1}-x_{t+1}^{j}),\; j=1,\dots,N.\]
and $J_{t+1}(x_{t+1})$ is approximated as \[
\hat{J}_{t+1}(x_{t+1})=\max_{j=1,\dots,N}\left\{ J_{t+1}(x_{t+1}^{j})+\nabla J_{t+1}(x_{t+1}^{j})^{T}(x_{t+1}-x_{t+1}^{j})\right\} .\]

We now show how to efficiently obtain $J_{t}(x_{t})$ and $\nabla J_{t}(x_{t})$
to approximate $J_{t}(x_{t})$.

Select a sample $\{w_{t}^{k}\}_{k=1}^{K}$ from the theoretical or
the empirical distribution of $w_{t}$, and rewrite (\ref{eq:Bellman})
as \begin{align}
J_{t}(x_{t})= & \inf_{u_{t}}\sum_{k=1}^{K}P(w_{t}^{k})[c(x_{t},u_{t},w_{t}^{k})+J_{t+1}(x_{t+1}(x_{t},u_{t},w_{t}^{k}))]\nonumber \\
\text{s.t.}\quad & g_{t}(x_{t},u_{t})\le0.\label{eq:Bellman_sample}\end{align}
If $w_{t}$ is discrete and the set of it's possible values is small,
we can select all of its possible values as $\{w_{t}^{k}\}_{k=1}^{K}$;
otherwise $\{w_{t}^{k}\}_{k=1}^{K}$ will have to be a sample from
the distribution. The error introduced by this sampling will not be
considered in this paper, and (\ref{eq:Bellman_sample}) will be treated
as the definition of $J_{t}(x_{t})$.

To obtain $J_{t}$ and $\nabla J_{t}$ at a point $x_{t}$, we need
to solve an optimization problem. Substitute $J_{t+1}$ with its approximation
$\hat{J}_{t+1}$. To simplify the evaluation of $\hat{J}_{t+1}(x_{t+1}(x_{t},u_{t},w_{t}^{k}))$,
we replace it with a decision variable $J_{t+1}^{k}$, and add the
supporting hyperplane constraints (linear):\begin{align*}
\min_{u_{t},\{J_{t+1}^{k}\}} & \sum_{k=1}^{K}P(w_{t}^{k})[c(x_{t},u_{t},w_{t}^{k})+J_{t+1}^{k}]\\
\text{s.t.}\quad & g_{t}(x_{t},u_{t})\le0,\\
 & J_{t+1}^{k}\ge J_{t+1}(x_{t+1}^{j})+\nabla J_{t+1}(x_{t+1}^{j})^{T}(x_{t+1}(x_{t},u_{t},w_{t}^{k})-x_{t+1}^{j}),\;\text{for all }j,k.\end{align*}
The above formulation allows us to get $J_{t}(x_{t})$, but to obtain
$\nabla J_{t}(x_{t})$, we need to change it a bit. We introduce decision
variable $s_{t}$ as a dummy of $x_{t}$:\begin{align}
\min_{s_{t},u_{t},\{J_{t+1}^{k}\}} & \sum_{k=1}^{K}P(w_{t}^{k})[c(s_{t},u_{t},w_{t}^{k})+J_{t+1}^{k}]\nonumber \\
\text{s.t.}\quad & g_{t}(s_{t},u_{t})\le0,\nonumber \\
 & J_{t+1}^{k}\ge J_{t+1}(x_{t+1}^{j})+\nabla J_{t+1}(x_{t+1}^{j})^{T}(x_{t+1}(s_{t},u_{t},w_{t}^{k})-x_{t+1}^{j}),\;\text{for all }j,k,\nonumber \\
 & s_{t}=x_{t}.\label{eq:Bellman_dummy}\end{align}
Solving (\ref{eq:Bellman_dummy}) not only gives us $J_{t}(x_{t})$,
which is the optimal objective value, but our assumptions also guarantee
that there exists a Lagrange multiplier vector $\lambda$ associated
to the constraint $s_{t}=x_{t}$, and that $-\lambda$ (or $\lambda$,
depending on how you write the Lagrangian function) is a subgradient
of $J_{t}()$ at $x_{t}$. Thus we can obtain $\nabla J_{t}(x_{t})$
for free by looking at the Lagrange multiplier of the constraint $s_{t}=x_{t}$.

Note that approximating a function with supporting hyperplanes preserves
convexity, thus (\ref{eq:Bellman_dummy}) is a convex program, and
any local minimum is also a global minimum.

\subsection*{Adaptive Convex Enveloping}

We just showed how to efficiently obtain $J_{t}$ and $\nabla J_{t}$
to build supporting hyperplanes to approximate $J_{t}(x_{t})$, but
we haven't discussed where to build these hyperplanes. From the geometric
point of view, the flatter the function in a region, the less supporting
hyperplanes needed, vice versa. In principle, we want to be as economical
as possible, because each supporting hyperplane will add a number
of linear constraints to (\ref{eq:Bellman_dummy}) in stage $t-1$.
Without prior knowledge of $J_{t}(x_{t})$, presetting the investigation
points $\{x_{t}^{j}\}$ could be wasteful in supporting hyperplanes
where the function is flat, and insufficient where the function is
very curved.

The approach used by A.C.E. --- the reason why it got the name ---
is to learn the shape of $J_{t}(x_{t})$ on the way, adding supporting
hyperplanes only where necessary.

\subsubsection*{Error Control}

The supporting hyperplanes give a lower bound of $J_{t}(x_{t})$ at
any $x_{t}\in X_{t}$, and the convexity of $J_{t}(x_{t})$ can give
an upper bound --- together they provide an error bound.

\begin{figure}
\begin{centering}
\includegraphics[width=0.8\textwidth]{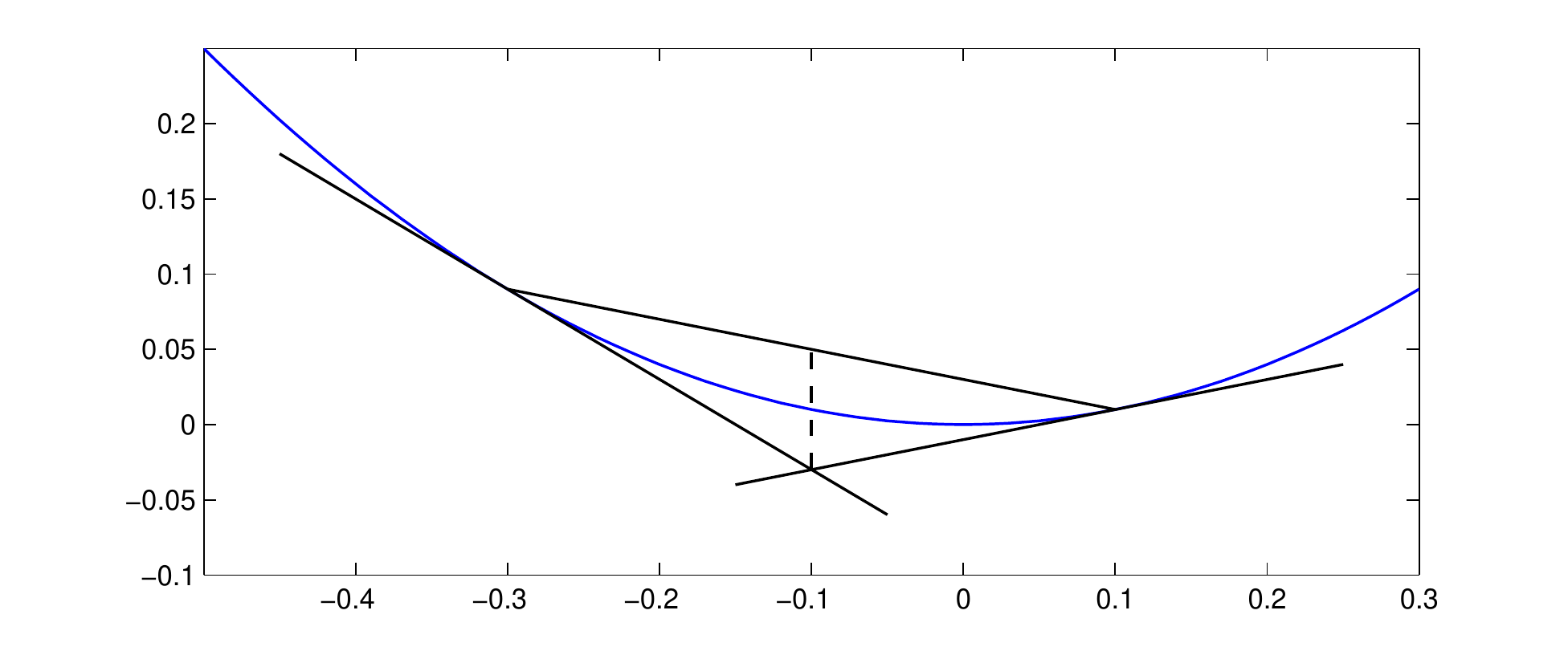}
\par\end{centering}

\caption{upper, lower bounds and max error}

\label{Flo:error_bound}
\end{figure}

Suppose that our state variable $x_{t}$ is 1-dimensional and look
at Figure \ref{Flo:error_bound}. We have two tangents added at $-0.3$
and $0.1$, respectively. By convexity, for $x_{t}\in[-0.3,0.1]$,
$J_{t}(x_{t})$ must be above the two tangents, and below the segment
connecting the two points of tangency. Thus if we use the max of the
tangents as the approximation of $J_{t}(x_{t})$, the maximum potential
error at $x_{t}$ is the vertical distance from the max of the tangents
to the segment connecting the two points of tangency. And the error
of any point in the region $[-0.3,0.1]$ is bounded by the maximum
potential error at the intersection of the two tangents.

In general, suppose that the domain $X_{t}$ is a $p$-dimensional
body, and let $\{x_{t}^{j}\}_{j=1}^{p+1}$ be $p+1$ points whose
convex hull $H_{conv}(\{x_{t}^{j}\}_{j=1}^{p+1})$ is also $p$-dimensional.
If we add supporting hyperplanes at $\{x_{t}^{j}\}$, $j=1,\dots,N$
($N\ge p+1$), and use the max of these hyperplanes as the approximation
of $J_{t}(x_{t})$, then for any $x_{t}\in X_{t}\cap H_{conv}(\{x_{t}^{j}\}_{j=1}^{p+1})$,
the maximum potential error is the vertical distance from the max
of the supporting hyperplanes to the hyperplane defined by the points
$\{(x_{t}^{j},J_{t}(x_{t}^{j})\}_{j=1}^{p+1}$. The maximum potential
error for the region $X_{t}\cap H_{conv}(\{x_{t}^{j}\}_{j=1}^{p+1})$
can be found by expressing the point $x_{t}$ as a convex combination
of $\{x_{t}^{j}\}_{j=1}^{p+1}$ and solving the following maximization
problem: \begin{align}
\max_{x_{t},\alpha,J}\; & \sum_{j=1}^{p+1}\alpha_{j}J_{t}(x_{t}^{j})-J\nonumber \\
\text{s.t.}\; & x_{t}=\sum_{j=1}^{p+1}\alpha_{j}x_{t}^{j},\nonumber \\
 & x_{t}\in X_{t},\nonumber \\
 & \alpha_{j}\ge0,\;\text{for }j=1,\dots,p+1,\nonumber \\
 & \sum_{j=1}^{p+1}\alpha_{j}=1,\nonumber \\
 & J\ge J_{t}(x_{t}^{j})+\nabla J_{t}(x_{t}^{j})^{T}(x-x_{t}^{j}),\;\text{for }j=1,\dots,N.\label{eq:max_err_seeker}\end{align}
There are two facts that may counter one's intuition when $p>1$.
First, the intersection of the supporting hyperplanes at $\{x_{t}^{j}\}_{j=1}^{p+1}$
may not be in $H_{conv}(\{x_{t}^{j}\}_{j=1}^{p+1})$; it may not be
in $X_{t}$, too, even when $H_{conv}(\{x_{t}^{j}\}_{j=1}^{p+1})\subset X_{t}$.
Second, a supporting hyperplane at $x_{t}^{j}$, $j>p+1$, can be
an active lower bound, even when $x_{t}^{j}\notin H_{conv}(\{x_{t}^{j}\}_{j=1}^{p+1})$.
These are the reasons why we need to solve (\ref{eq:max_err_seeker})
to find the potentially worst point $x_{t}$, and why we use all the
supporting hyperplanes as constraints.

\subsubsection*{Recursive Partitioning}

\begin{figure}
\begin{centering}
\includegraphics[width=0.45\textwidth]{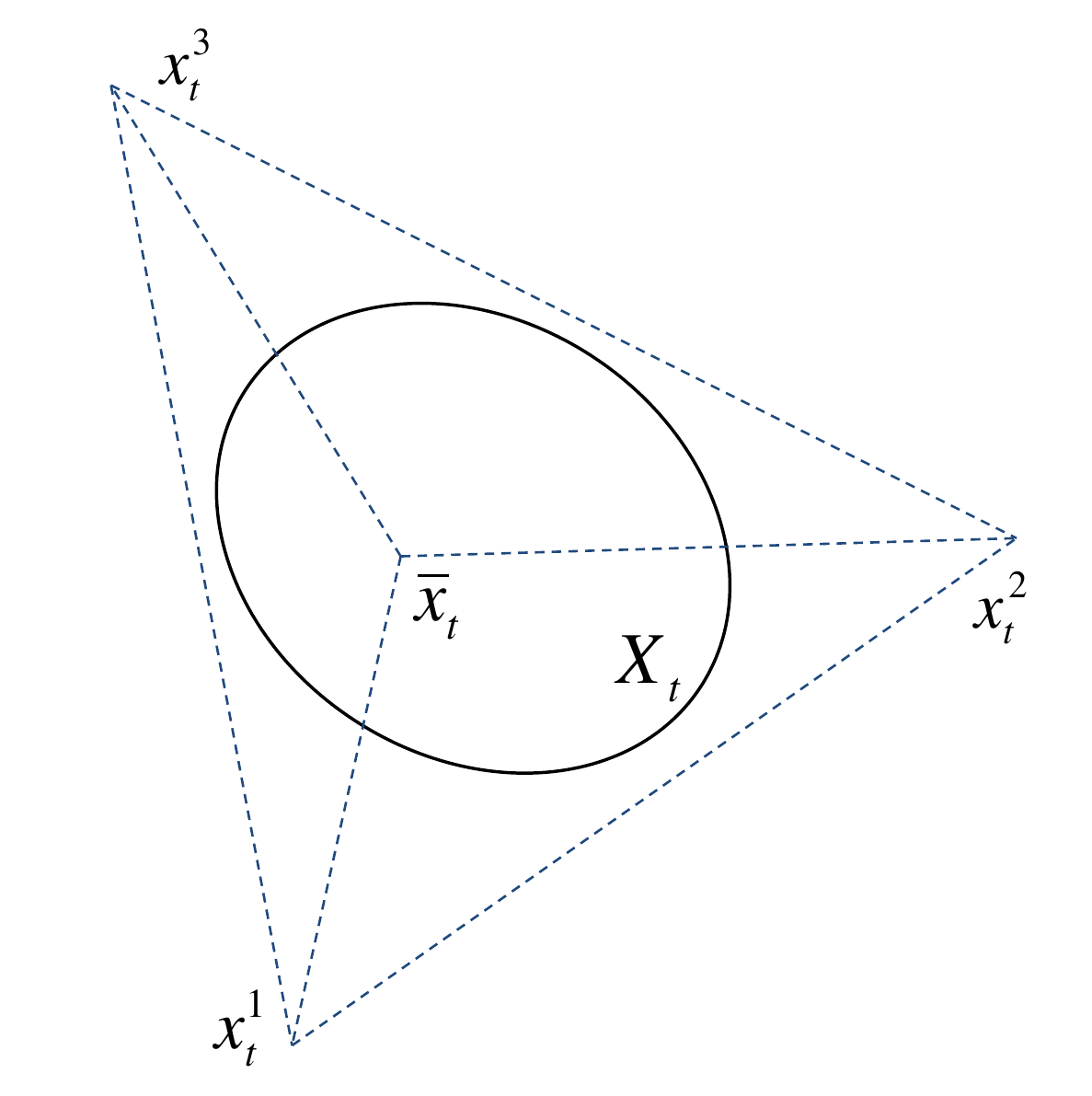}
\par\end{centering}

\caption{Partitioning at the potentially worst point}

\label{Flo:partitioning}
\end{figure}

Let the optimal solution of (\ref{eq:max_err_seeker}) be $(\bar{x}_{t},\bar{\alpha},\bar{J})$.
The point $\bar{x}_{t}$ is where the maximum potential error could
occur. If this potential error is larger than the tolerance, we can
add a supporting hyperplane at $\bar{x}_{t}$ to reduce the error,
and separate $H_{conv}(\{x_{t}^{j}\}_{j=1}^{p+1})$ at $\bar{x}_{t}$
into $p+1$ smaller convex hulls, each with $\bar{x}_{t}$ and $p$
points from $\{x_{t}^{j}\}_{j=1}^{p+1}$ as vertices (See Figure \ref{Flo:partitioning}).
Then we face $p+1$ subregions whose maximum potential errors can
be found by solving (\ref{eq:max_err_seeker}) again. Note that if
$\bar{x}$ is on a facet or an edge of $H_{conv}(\{x_{t}^{j}\}_{j=1}^{p+1})$,
some of the subregions will be less than $p$-dimensional. We ignore
any subregion that is less than $p$-dimensional, as it is covered
by the other $p$-dimensional subregions. We can tell if a subregion
is $p$-dimensional by looking at $\bar{\alpha}$.

To initiate the recursive partitioning algorithm, we need to have
$p+1$ initial vertices $\{x_{t}^{j}\}_{j=1}^{p+1}$ that satisfy
two conditions: 1, their convex hull $H_{conv}(\{x_{t}^{j}\}_{j=1}^{p+1})$
contains $X_{t}$; 2, the feasible action sets $\{u_{t}\mid g_{t}(x_{t}^{j},u_{t})\le0\}$,
$j=1,\dots,p+1$, are nonempty. The first condition is for controlling
the error at any point of $X_{t}$, while the second condition is
needed so that we can build supporting hyperplanes at these vertices.
These points are usually not hard to find. For example, if the states
have constraints $0\le x_{t,i}\le1$, $i=1,\dots,p$, one may want
to check the intersections of the hyperplanes $x_{i}=0$, $i=1,\dots,p$,
and $\sum_{i=1}^{p}x_{i}=p$ to see if they satisfy the second condition.
However, it's not guaranteed that these points always exist or are
always easy to find. If one really cannot find these points and has
no way around it, then one may want to make some compromise and choose
$\{x_{t}^{j}\}_{j=1}^{p+1}$ that satisfy the second condition and
their convex hull covers as much of $X_{t}$ as possible.

Once we have the initial convex hull to start the algorithm, we separate
it into sub-hulls if its maximum potential error exceeds the tolerance.
Repeat this procedure recursively to all the sub-hulls, until all
of their maximum potential errors are less than or equal to the tolerance.

\begin{algorithm}[H]
\begin{lyxcode}
Loop~t~from~T-1~to~1
\begin{lyxcode}
Read~$J_{t+1}$~supporting~hyperplane~information~from~file;

Find~initial~vertices~$\{x_{t}^{j}\}_{j=1}^{p+1}$;

Add~supporting~hyperplanes~at~$\{x_{t}^{j}\}_{j=1}^{p+1}$~by~solving~(\ref{eq:Bellman_dummy});

Form~the~first~item~of~the~section~list~with~$\{x_{t}^{j}\}_{j=1}^{p+1}$;

currentMaxError~=~tolerance~+~1;

While~(currentMaxError~>~tolerance)
\begin{lyxcode}
If~(list~size~>~Budget)
\begin{lyxcode}
Print~{}``Budget~exceeded.'';

Break;
\end{lyxcode}
End

currentMaxError~=~0;

iterator~=~beginning~of~the~list;

While~(iterator~!=~end~of~the~list)
\begin{lyxcode}
Solve~(\ref{eq:max_err_seeker});

currentMaxError~=~max(currentMaxError,~optimal~value~of~(\ref{eq:max_err_seeker}));

If~(optimal~value~of~(\ref{eq:max_err_seeker})~<~tolerance)
\begin{lyxcode}
iterator++;
\end{lyxcode}
Else
\begin{lyxcode}
Solve~(\ref{eq:Bellman_dummy});

Add~a~new~supporting~hyperplane~at~$\bar{x}_{t}$;

Separate~the~current~section~at~$\bar{x}_{t}$~to~subsections;

tmp\_iterator~=~iterator;

iterator++;

Replace~tmp\_iterator~with~new~subsections;
\end{lyxcode}
End
\end{lyxcode}
End
\end{lyxcode}
End

Write~$J_{t}$~supporting~hyperplane~information~to~file;

Release~memory;
\end{lyxcode}
End
\end{lyxcode}
\caption{Recursive Partitioning}

\label{Flo:algorithm_partition}
\end{algorithm}

Algorithm \ref{Flo:algorithm_partition} shows how to solve the dynamic
program with recursive partitioning. In the algorithm, we call a convex
hull a {}``section'', and we use a list of sections. The beginning
of a list is the list's first item, while the end of a list is a position
beyond the last item. An iterator of a list is like a pointer that
points to an item of the list. The {}``\texttt{\textsc{++}}'' operator
points the iterator to the next item, but if iterator is the last
item of the list, {}``\texttt{iterator++}'' will point it to the
end of the list, i.e., a position beyond the last item.

Note that the error we talk about here is the error of $\hat{J}_{t}$
relative to $\hat{J}_{t+1}$. The absolute error of $\hat{J}_{t}$
is the error relative to $\hat{J}_{t+1}$ plus the absolute error
of $\hat{J}_{t+1}$.

\subsubsection*{Approximating by Importance}

Theoretically, we can approximate a function to any arbitrary high
precision by reducing the tolerance. In practice, however, doing so
to a high dimensional function is prohibitive. Meanwhile, pursuing
a high precision at every point in $X_{t}$ could be wasteful, since
in many applications the optimal policy $\pi$ tends to guide the
process to visit frequently only a small portion of the domain. Therefore
it is reasonable to only focus on these more important small portions.

We don't know where these small portions are beforehand, however.
What we could do is to use recursive partitioning to obtain affordable
and relatively good approximations of the cost-to-go functions, and
use simulation to see where the policy guides the process. The following
is one possible way to update the approximations by importance.
\begin{lyxcode}
Observe~state~path~$\{x_{t}\}_{t=1}^{T}$~under~the~currect~policy;

Loop~t~from~T-1~to~1
\begin{lyxcode}
Find~the~section~that~contains~$x_{t}$;

If~(maximum~potential~error~>~tolerance)
\begin{lyxcode}
Add~a~supporting~hyperplane~at~$x_{t}$;

Separate~the~current~section~at~$x_{t}$;
\end{lyxcode}
End
\end{lyxcode}
End
\end{lyxcode}
This refining process can run as long as needed. Note that this recursive
partitioning \& approximating by importance approach bypasses the
exploration vs.~exploitation dilemma, which troubles many of today's
popular ADP methods.

\subsection*{Example}

Here is a demonstration that applies A.C.E.~to the well-known inventory
control problem. Let the purchasing cost, penalty and holding cost
be 2.0, 4.0 and 0.2, respectively. Let the demand be uniform on $[0,10]$,
and let the sample be 0.0, 0.1,$\dots$,9.9. The cost-to-go function
is known to be convex, which also could be confirmed by checking the
convexity conditions given at the beginning of this paper. The optimal
policy is well known, that is, to order and increase the inventory
up to a level $S$, but the value of $S$ is not known, and need to
be computed numerically.

The following are the approximations given by A.C.E.~for $x\in[0,15]$,
$t$ from $T-1$ to $T-10$, at tolerance = 0.1. $J_{T}$ is assumed
to be zero.

\begin{figure}
\begin{centering}
\includegraphics[width=0.6\textwidth]{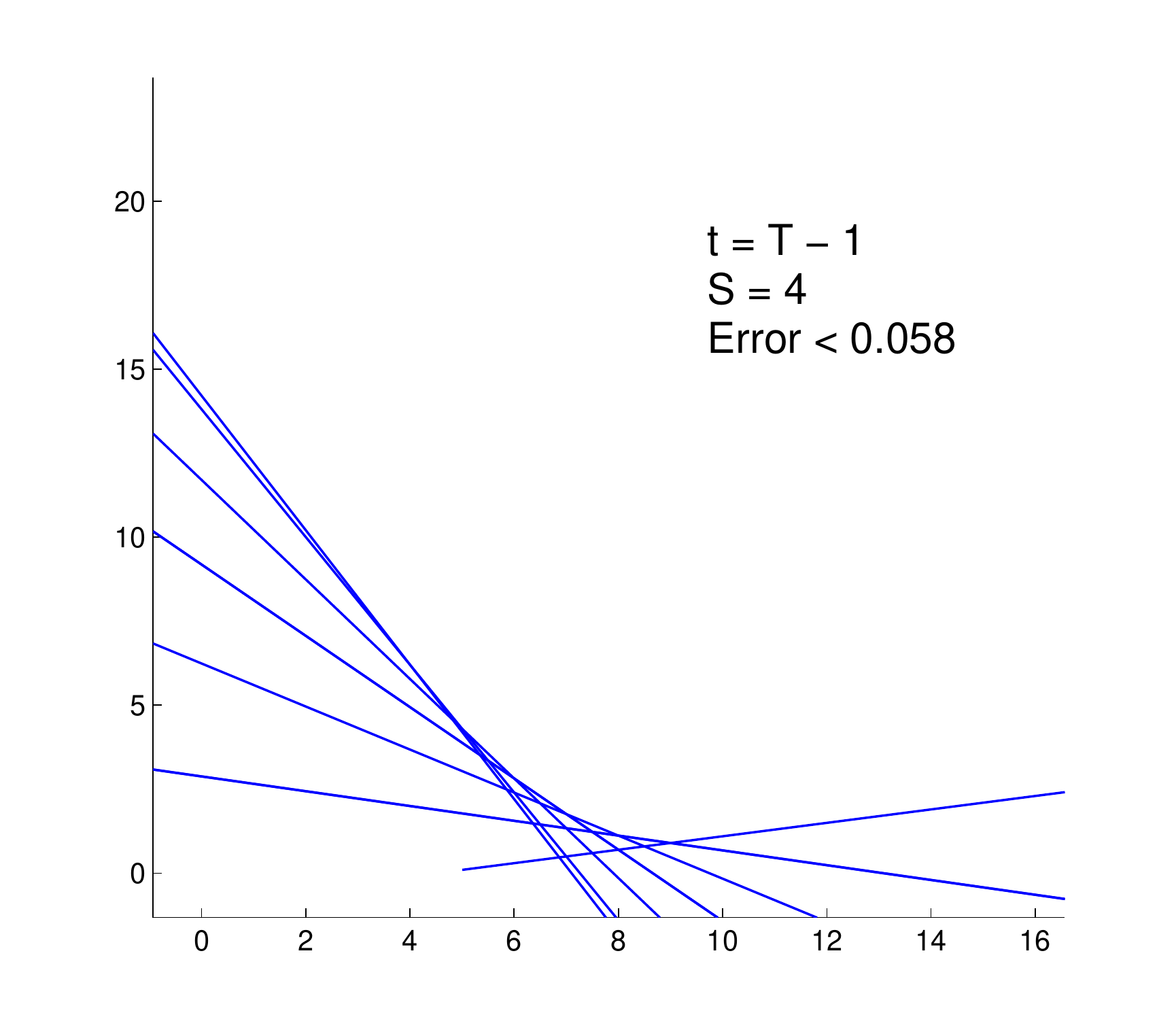}
\par\end{centering}

\caption{A.C.E. approximation of $J_{T-1}$}

\label{Flo:J_T-1}
\end{figure}

The approximation of $J_{T-1}$ is shown by Figure \ref{Flo:J_T-1}.
The optimal policy is to order up to $S=4$. For the region $x\le4$
and $x\ge10$, we know $J_{T-1}(x)$ should be linear, and we see
A.C.E.~didn't waste supporting hyperplanes (tangents) there. For
the region $4\le x\le10$, it can be shown that $J_{T-1}(x)$ is quadratic
under a demand from uniform distribution, and the approximation reflects
the shape closely. Overall, the approximation has an error less than
0.058 at any point $x\in[0,15]$.

\begin{figure}
\begin{centering}
\includegraphics[width=0.45\textwidth]{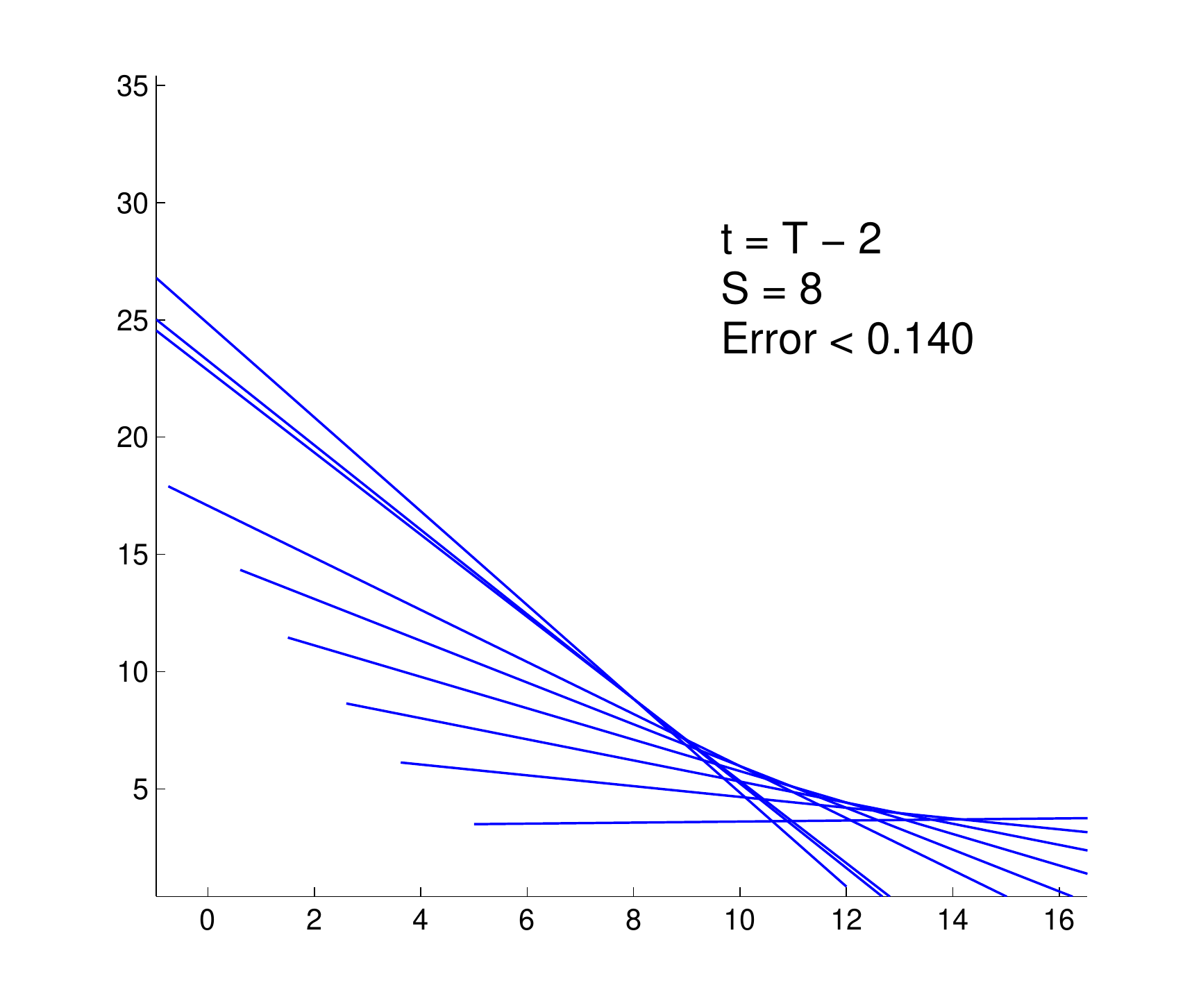}\includegraphics[width=0.45\textwidth]{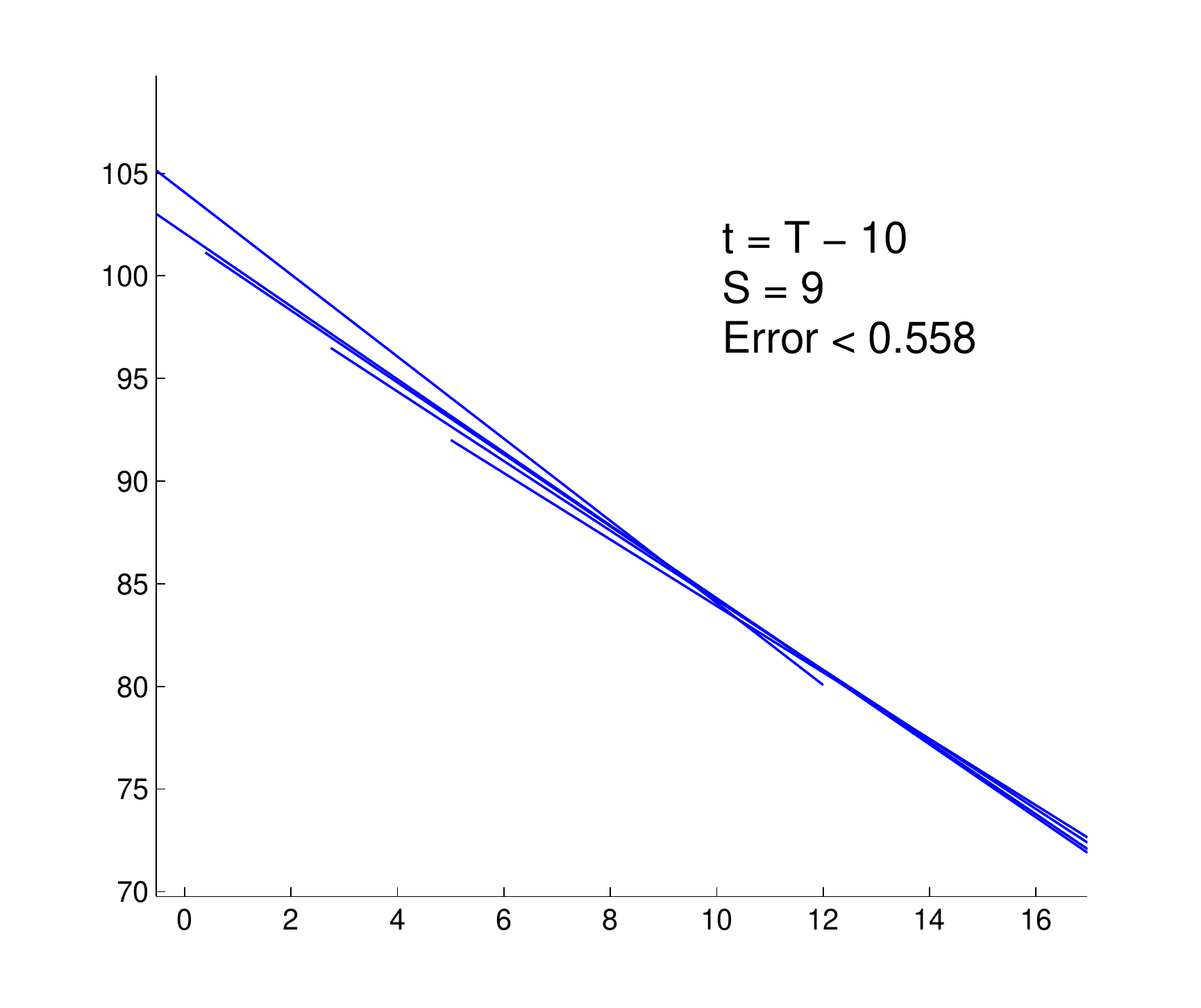}
\par\end{centering}

\caption{A.C.E. approximation of $J_{T-2}$ and $J_{T-10}$}

\label{Flo:J_T-2_T-10}
\end{figure}

Figure \ref{Flo:J_T-2_T-10} shows the approximations of $J_{T-2}$
and $J_{T-10}$. The optimal policy at stage $T-2$ is to order up
to $S=8$. Starting (backwards) from stage $T-3$, the optimal policy
becomes fixed, and we always order up to $S=9$. The error accumulates.
At stage $T-10$, the error bound for a point $x\in[0,15]$ is 0.558.
But, the cost-to-go accumulates, too, and compared to the scale of
this target, the error is pretty small.

The whole process took 0.87 second on my computer. Since this is a
simple 1-dimensional example, approximating by importance is not needed.

\subsection*{Features of A.C.E.}
\begin{itemize}
\item Solves convex DP with multivariate continuous states and actions;
\item It is a standardized, general purpose method (no parameter tunning,
basis function choosing, kernel designing, etc.);
\item No assumption on the form of the cost-to-go functions --- the form
is learned on the way;
\item Implementation of mathematical programming allows large number of
action variables;
\item The supporting hyperplane approximation preserves convexity of the
cost-to-go functions, enabling reliable optimization;
\item The supporting hyperplane approximation also {}``ignores'' the state
variables that do not add to the nonlinearity of the target function,
e.g., approximating $f(x_{1},x_{2},x_{3})=x_{1}^{2}+x_{3}$ is no
more difficult than approximating $f(x)=x^{2}$, although the dimension
is higher;
\item The computation is relatively light since we obtain the subgradients
for free from the Lagrange multipliers.
\item The supporting hyperplanes are constructed economically by need and
importance, further reducing the computation;
\item There are also problem dependent techniques that could greatly speed
up the computation. For instance, if $c_{t}$ and $g_{t}$ are piecewise
linear, we can use dual simplex updates to obtain the supporting hyperplanes,
instead of solving (\ref{eq:Bellman_dummy}) from scratch, because
(\ref{eq:Bellman_dummy}) is now a linear program and the state $x_{t}$
is a RHS;
\item We know an upper bound (not the big-$O$ notation) of the error between
the approximated and the true cost-to-go functions at any state, thus
we know how the policy performs;
\item There is no exploration vs.~exploitation dilemma;
\item Compared to stochastic programming, it stands out when solving problems
with long planning horizons, as the CPU time grows linearly in the
number of stages, and memory consumption remains the same; plus, it
is dynamic.
\end{itemize}

\end{document}